\documentclass[12pt]{article}
\usepackage{amsfonts}
\usepackage{amssymb}
\usepackage{amsthm}
\usepackage{newlfont}
\usepackage{amsmath,amsthm}

\makeatletter
 \@addtoreset{equation}{section}
\makeatother
\newtheorem{thm}{Theorem}[section]
\newtheorem{lem}{Lemma}[section]
\newtheorem{cor}{Corollary}[section]
\newtheorem{rema}{Remark}[section]
\newtheorem{prop}{Proposition}[section]
\newtheorem{Def}{Definition}[section]

\begin{document}

\title{The structure of stable  constant mean curvature hypersufaces }
\author{{Xu Cheng}\thanks{Supported by CNPq of Brazil}, Leung-fu Cheung  and Detang
Zhou\thanks{Supported by CAPES and CNPq  of Brazil.}} \maketitle

\begin{abstract}
We study the global behavior of (weakly) stable constant mean
curvature hypersurfaces in general Riemannian manifolds. By using
harmonic function theory, we prove some one-end theorems which are
new even for constant mean curvature hypersurfaces  in space forms.
 In particular,  a complete oriented weakly stable minimal hypersurface in $\mathbb{R}^{n+1}, n\geq 3,$
  must have only one end. Any complete noncompact weakly stable CMC $H$-hypersurface in
the hyperbolic space $\mathbb{H}^{n+1}, n=3,4,$ with
$H^2\geq\frac{10}{9}, \frac{7}{4},$ respectively, has only one end.

\end{abstract}

\section*{ 0 Introduction}\label{introduction}

The classical Bernstein theorem states that a  minimal entire graph
in $\mathbb{R}^3$ must be planar. This theorem was later generalized
to higher dimensions (dimension of the ambient Euclidean space
$\mathbb{R}^{n+1}$ is no more than $8$) by Fleming\cite{fl},
Almgren\cite{a}, De Giorgi\cite{dg}, and Simons\cite{s}.  In
$\mathbb{R}^{n+1}$,  $n\geq 8$, the examples of nonlinear entire
graphs are given by Bombieri, de Giorgi and Giusti \cite{bdgg}.
Because of the stability of minimal entire graphs, one is naturally
led to the generalization of the classical Bernstein theorem to the
question of asking whether all stable minimal hypersurfaces in
$R^{n+1}$ are hyperplanes when $n \leq 7$.  In the case when $n=2$,
this problem was solved independently by do Carmo and Peng
\cite{dcp}; and Fischer-Colbrie and Schoen \cite{fs}.  For higher
dimension, this problem is still open. On the other hand, there are
some results about the structure of stable minimal hypersurfaces in
all $\mathbb{R}^{n+1}$. For instance, H. Cao, Y. Shen and S. Zhu
\cite{csz} proved that a complete stable minimal hypersurface in
$\mathbb{R}^{n+1}, n\geq 3,$ must have only one end.

If the ambient manifold is not the Euclidean space, Fischer-Colbrie
and Schoen \cite{fs} gave a classification for complete oriented
stable minimal surfaces in a complete oriented $3$-manifold of
nonnegative scalar curvature. Recently, Li and Wang \cite{lw1}
showed that a complete noncompact properly immersed stable minimal
hypersurface in a complete manifold of nonnegative sectional
curvature must either have only one end  or be totally geodesic and
a product of a compact manifold with nonnegative sectional curvature
and $\mathbb{R}$.

In this paper we study  hypersurfaces with constant mean curvature
$H$. Let us now fix terminologies and notations so as to our
theorems. In the sequel we will abbreviate constant mean curvature
hypersurfces by calling them CMC $H-$hypersurfaces and will allow
$H$ to vanish (hence the need of putting $H$ here). Instead of the
usual stability, we will consider a weaken form of stability, which
is in fact the natural one for CMC $H-$hypersurfaces in case $H \ne
0$. Intuitively, a CMC hypersurface is weakly stable if the second
variations are nonnegative for all compactly supported
enclosed-volume-preserving variations (see Definition \ref{def1} and
Remark \ref{strongweak}). This concept of weakly stable CMC
hypersurfaces was introduced by Barbosa, do Carmo and Eschenburg
 \cite{bdce}, to accounts for the fact that spheres are stable (see  \cite{bdce}). This weak stability
comes naturally from the phenomenon of soap bubbles and is related
to isoperimetric problems.  In \cite{ds}, da Silveira studied
complete noncompact weakly stable CMC surfaces  in $\mathbb{R}^3$ or
the hyperbolic space $\mathbb{H}^3$. He proved that complete weakly
stable  CMC surfaces in $\mathbb{R}^3$ are planes and hence
generalized the corresponding result of do Carmo and Peng
\cite{dcp}, Fischer-Colbrie and Schoen \cite{fs}. For $\mathbb{H}^3$
he shows that  only horospheres can occur when constant mean
curvature $|H|\ge 1$. For higher dimensions,  very little is known
about complete noncompact weakly stable CMC hypersurfaces.

 In this paper, we  study the global behavior of weakly
stable CMC hypersurface (including minimal case). First, we obtain
\begin{thm}\label{sectcurvature1} (Th.\ref{sectcurvature}) Let $N^{n+1}, n\geq 5,$
be a complete Riemannian manifold and $M$ be a complete noncompact
weakly stable immersed  CMC $H$-hypersurface in $N$. If one of the
following cases  occurs,

(1) when $n=5$, the sectional curvature of $N$ is nonnegative and
$H\neq 0$;

(2) when $n\geq 6$, the sectional curvature $\tilde{K}$ of $N$
satisfies $\tilde{K}\geq \tau>0$ and $H^2\leq
\frac{4(2n-1)}{n^2(n-5)}\tau,$ for some number $\tau>o$;

(3) when $n\geq 6$, the sectional curvature and the Ricci curvature
of $N$ satisfy $\tilde{K}\geq 0$, $\tilde{\mathrm{Ric}}\geq \tau
>0$, for some number $\tau
>0$, and $H=0$,

 then  $M$ has only one end.

\end{thm}

The reason for the restriction on dimensions of  CMC hypersurfaces
in the above theorem is that there are some nonexistence results
(see the proof of this theorem for detail). Theorem
\ref{sectcurvature1} has the following examples: complete noncompact
weakly stable CMC $H$-hypersurfaces in    the standard sphere
$\mathbb{S}^6$ with $H\neq 0$;  or in the standard sphere
$\mathbb{S}^{n+1}, n\geq 6$ with $H^2\leq \frac{4(2n-1)}{n^2(n-5)}$.

Actually, Theorem \ref{sectcurvature1} is a special case of the more
general Theorem \ref{one}, which also implies that

\begin{thm}\label{oneHn1}(Cor.\ref{oneHn}) Any complete noncompact weakly stable CMC $H$-hypersurface in
the hyperbolic space $\mathbb{H}^{n+1}, n=3,4,$ with
$H^2\geq\frac{10}{9}, \frac{7}{4},$ respectively, has only one end.
\end{thm}

Next we  consider complete weakly stable  minimal hypersurfaces in
$\mathbb{R}^{n+1},$ $n\geq 3,$  and generalize the results of Cao,
Shen and Zhu as follows:
\begin{thm}\label{oneRn1}(Th.\ref{oneRRn}) A complete oriented weakly stable minimal
hypersurface in $\mathbb{R}^{n+1}, n\geq 3,$ must have only one end.
\end{thm}

With this theorem, we obtain
\begin{cor} \label{plane1}(Cor.\ref{plane}) A complete oriented weakly stable immersed
minimal hypersurface in $\mathbb{R}^{n+1}, n\geq 3,$ with finite
total  curvature (i.e., $\int_M|A|^n<\infty$) is a hyperplane.
\end{cor}

Finally, we study the structure of weakly stable CMC hypersurfaces
according to the parabolicity or nonparabolicity of $M$. We obtain
the following results:

\begin{thm}\label{thmpara1}(Th.\ref{thmpara}) Let $N$ be a complete  manifold of bounded geometry and $M$ be
a complete noncompact  weakly stable CMC $H$-hypersurface immersed
in $N$. If the sectional curvature of $N$ is bounded from below by
$-H^2$ and  $M$ is parabolic, then it is totally umbilic and has
nonnegative sectional curvature. Furthermore,  either

(1) $M$ has only one end; or

(2) $M=\mathbb{R}\times P$ with the product metric, where $P$ is a
compact manifold of nonnegative sectional curvature.

\end{thm}

\begin{thm}\label{thmnonpara1}(Th.\ref{thmnonpara}) Let $N$ be a complete Riemannian manifold and $M$ be
a complete noncompact  weakly stable CMC $H$-hypersurface immersed
in $N$. If $M$ is nonparabolic, and
$$\tilde{\mathrm{Ric}}(\nu)+\tilde{\mathrm{Ric}}(X)-\tilde{\mathrm{K}}(X,\nu)\geq \frac{n^2(n-5)}{4}H^2, \forall X\in T_pM, |X|=1, p\in M,$$

then it has only one nonparabolic end, where $\tilde{\mathrm{K}}$
and $\tilde{\mathrm{Ric}}$  denote the sectional and Ricci
curvatures of $N$, respectively; $\nu$ denotes the unit normal
vector field of $M$.
\end{thm}

In  some of recent works, the structure of stable (i.e., strongly
stable) minimal hypersurfaces was studied by means of harmonic
function theory (see \cite{csz}, \cite{lw}, \cite{lw1}). The same
approach can be used in the study of weakly stable CMC
hypersurfaces. However, a significant difference between weakly
stable and strongly stable cases lies in the choice of test
functions. When one deals with weak stability, the test functions
$f$ must satisfy $\int_M f=0$. In this paper, we successfully
construct the required test functions by using the properties of
harmonic functions (Theorem \ref{thm3}  and Proposition \ref{sch}).
Combining our construction and the approach in \cite{lw},
\cite{lw1}, we are  able to discuss the global behavior of weakly
stable CMC hypersurfaces.  In Theorem \ref{thm3}, we obtain the
nonexistence of nonconstant bounded harmonic functions with finite
Dirichlet integral on weakly stable CMC hypersurfaces. This theorem
enable us to study the uniqueness of ends. In Proposition \ref{sch},
we discuss a property of Schr\"odinger operator on parabolic
manifolds
 which can be applied to study weakly stable
CMC hypersurfaces with parabolicity. Besides, different from minimal
hypersurfaces, CMC hypersurfaces with $H\neq 0$ have  the curvature
estimate depending on $H$, which causes dimension restriction in the
results.

The rest of this paper is organized as follows: in Section 1 we give
some definitions and facts as preliminaries; in Section 2, we first
discuss volume growth of the ends of complete noncompact
hypersurfaces with mean curvature vector field bounded in norm, and
then  study nonparabolicity  of the ends of CMC hypersurfaces with
stability; in Section 3, we use  harmonic functions to study the
uniqueness of ends of complete noncompact weakly stable CMC
hypersurfaces; in Section 4, we give a property of Schr\"odinger
operator on parabolic manifolds; in the last section (Section 5), we
discuss the structure of complete noncompact weakly stable CMC
hypersurfaces.

The results on minimal case in this paper has been announced in
\cite{ccz}.

\bigskip
\noindent\textbf{Acknowledgements.} One part of this work was done
while the third author was visiting the Department of  Mathematics,
University of California, Irvine. He wishes to thank the department
 for  hospitality. The authors would like to thank Peter Li for
 some conversations.

\section{Preliminaries}\label{pre}

We recall some definitions and facts in this section.

Let $N^{n+1}$ be an oriented $(n+1)$-dimensional Riemannian manifold
and let $i:M^n\to N^{n+1}$ be an  isometric immersion of a connected
$n$-dimensional manifold $M$ with constant mean curvature $H$. We
assume $M$ is orientable. When $H$ is nonzero, the orientation is
automatic. Throughout this paper, $\tilde{\mathrm{K}},
\tilde{\mathrm{Ric}}$, $\mathrm{K}$, and $\mathrm{Ric}$  denote the
sectional, Ricci curvatures of $N$, the sectional, Ricci curvature
of $M$ respectively. $\nu$ denotes the unit normal vector field of
$M$. $|A|$ is the norm of the second fundamental form $A$. $B_p(R)$
will denote the  intrinsic geodesic ball in $M$ of radius $R$
centered at $p$. We have

\begin{Def} \label{def1} There are two cases.
In the case $H\neq 0$,
the immersion $i$ is called {\it stable} or {\it weakly stable} if

\begin{equation}\label{eqnstable}
\int_{M}\{|\nabla f|^{2}-\left( \tilde{\mathrm{Ric}}(\nu,\nu
)+|A|^{2}\right) f^{2}\}\geq0,
\end{equation}
for all compactly supported piecewise smooth  functions $%
f\colon M\rightarrow R$ satisfying
\[
\int_{M}f=0,
\]
where  $\nabla f$ is the gradient of $f$ in the induced metric of
$M$;

the immersion $i$ is called {\it strongly stable} if
(\ref{eqnstable}) holds for all compactly supported piecewise smooth
functions
$%
f\colon M\rightarrow R$.

In the case $H=0$ (minimal case),  the immersion $i$ is called {\it weakly stable} if (\ref{eqnstable}) holds for all compactly supported piecewise
 smooth functions satisfying
$%
f\colon M\rightarrow R$
\[
\int_{M}f=0;
\]

the immersion $i$ is called {\it stable} if (\ref{eqnstable}) holds
for all compactly supported piecewise smooth functions
$%
f\colon M\rightarrow R$.

\end{Def}
It is known, from the definition, that a weakly stable minimal
hypersurface has the index 0 or 1 (see \cite{ds}).  Obviously, a
strongly stable CMC hypersurface is weakly stable. But the converse
may not be true. For example, $\mathbb{S}^2\subset \mathbb{S}^3$ as
a totally geodesic embedding in the ordinary $3$-sphere is not
stable but weakly stable.

\begin{rema}\label{strongweak} In the current literatures, the terms of
stability on minimal  and constant mean curvature hypersurfaces are
different (perhaps a little confusing). A hypersurface with nonzero
constant mean curvature is called stable if it is weakly stable;
while a minimal hypersurface is called stable if it is strongly
stable in the above sense. In this paper, we deal with the weak
stability for both hypersurfaces. In order to avoid confusion and
conform to the notations of others, the notation of weak stability
is used without omission in this paper.

\end{rema}

For CMC $H$-hypersurfaces, it is convenient to introduce the
(traceless) tensor $\Phi:=A-HI$, where $I$ denotes the identity. A
straightforward computation gives $|\Phi|^2=|A|^2-nH^2$ and the
stability inequality (\ref{eqnstable}) becomes
\begin{equation}\label{eqnstabletr}
\int_{M}\{|\nabla f|^{2}-\left( \tilde{\mathrm{Ric}}(\nu,\nu
)+|\Phi|^{2}+nH^2\right) f^{2}\}\geq 0.
\end{equation}

In this paper, we will discuss the number of ends of hypersurfaces.
Now we give some related definitions.

\begin{Def} (cf. \cite{lt}, \cite{lw}) A manifold is said to be parabolic if it does not admit a positive Green's function. Conversely, a nonparabolic
manifold is one which admits a positive Green's function.

An end $E$ of $\Sigma$ is said to be  nonparabolic if it admits a
positive Green´s function with Neumann boundary condition on
$\partial E$. Otherwise, it is said to be parabolic.

\end{Def}

In order to estimate the number of ends of a weakly stable CMC
hypersurface,  we need the following theorem by Li and Tam.

\begin{thm}\label{endform}(\cite{lt}, see also\cite{lw} Theorem $1$) Let $M$ be a complete Riemannian manifold. Let $\mathcal{H}^0_D(M)$ be
the space of bounded harmonic functions with finite Dirichlet
integral. Then the number of nonparabolic ends of $M$ is bounded
from above by $\text{dim}\mathcal{H}^0_D(M)$.
\end{thm}

From Theorem \ref{endform}, we know that if every end of $M$ is
nonparabolic,  then  the number of its ends is no more than $\dim
\mathcal{H}^0_D(M)$.

\section{Nonparabolicity of ends}
In this section, we first discuss the volume growth of ends of
complete noncompact submanifolds in a Riemannian manifold $N$ of
bounded geometry ({\it a manifold $N$ is called bounded geometry if
its sectional curvatures  $\tilde{\mathrm{K}}\leq\sigma^2, \sigma>0$
and  its injectivity radius  $i_N(p)\geq i_0, i_0>0$}) and using it
to study the property of the nonparabolic ends of submanifolds.

Frensel \cite{fr} showed  that if $M$ is a complete noncompact
immersed submanifold in a
 manifold of bounded geometry with mean curvature vector field bounded in norm, then $M$ has infinite volume. Here, we  prove that even each end of $M$ has infinite volume.
\begin{lem}\label{volesti}(\cite{fr} Th.3) Let $N$ be an $m$-dimensional  manifold  and let $M$ be an $n$-dimensional complete noncompact  manifold. Let $x: M^n\rightarrow N^m$ be an isometric
immersion with mean curvature vector field bounded in norm. Assume
that $N$ has sectional curvature $\tilde{\mathrm{K}}\leq \sigma^2$,
where constant $\sigma>0$. Then
$$\mathrm{Vol}(B_p(R))\geq \sigma^{-n}\omega_n(\sin R\sigma)^ne^{-H_0R},$$
where $R\leq\min\{\frac{\pi}{2\sigma}, i_N(p)\}$ and $|H|\leq H_0$.
\end{lem}

We obtain
\begin{prop}\label{infvol}Let $N$ be an $m$-dimensional  manifold of bounded geometry  and let $M$ be an $n$-dimensional complete noncompact  manifold. Let $x: M\rightarrow N$ be an isometric
immersion with mean curvature vector field bounded in norm. Then
each end $E$ of $M$ has infinite volume. More exactly, the rate of
volume growth of $E$ is at least linear, i.e., for any $p\in E$,
\begin{equation}
\liminf_{R\rightarrow\infty}\frac{\mathrm{Vol}(B_p(R)\cap E)}{R}>0,
\end{equation}
where the limit is independent of the choice of $p$.

\end{prop}

\begin{proof} Assume that $E$ is  an end of $M$ with respect to  a compact set $D\subset M$ with smooth boundary $\partial D$.

{\it We claim that there exist some  $x\in E$ and a ray $\gamma$ in
$E$ emanating from $x$, i.e., $\gamma: [0,\infty)\rightarrow E$ is a
minimizing geodesic  satisfying $\gamma(0)=x$, and
$d(\gamma(s),\gamma(t))=|s-t|$, for all $s, t\geq 0$, where $\gamma$
has the arc length parameter.}

Now we prove the claim. Since $E$ is unbounded, there exists a
sequence of points $q_i\in E$ such that $d(q_i, D)\rightarrow
\infty$ when $i\rightarrow \infty$.  Since $D$ is compact, there
exist a sequence of points $p_i\in \partial D$ and  a sequence of
minimizing normalized  geodesic segments $\gamma_i|_{[0,s_i]}$ in
$M$ joining $p_i=\gamma_i(0)$ to $q_i=\gamma(s_i)$ respectively,
such that  $d(q_i,p_i)=d(q_i, D)$.

Each $\gamma_i$ has the following properties: 1) $p_i$ is the only
intersection of $D$ and $\gamma_i$ (otherwise, $d(q_i,p_i)\neq
d(q_i,D)$); 2) $\gamma_i\backslash \{p_i\}\subset E$ (since $E$ is a
connected component of $M\backslash D$); 3) $\gamma'_i(0)$ is
orthogonal to $D$ at $p_i$ (since $\gamma_i$ realizes the distance
$d(q_i,D)$).

Since the unit normal bundle of $D$ is compact, there exists a
subsequence of $(p_i,\gamma_i'(0))$, which is still denoted by
$(p_i,\gamma_i'(0))$, converging to a point $(p_0,\nu)$ in the unit
normal bundle, where $p_0\in D, \nu\in T_{p_0}M$ . Let
$\tilde{\gamma}|_{[0,+\infty)}$ be the normalized geodesic in $M$
emanating from $p_0$ with initial unit tangent vector $\nu$. By ODE
theory,  $\gamma_i$ converges to $\tilde{\gamma}$ uniformly on any
compact subset of $[0,+\infty)$. Moreover,  for any $s\in
[0,+\infty)$,  the segment $\tilde{\gamma}|_{[0,s]}$ realizes the
distance from $\tilde{\gamma}(s)$ to $D$. By the same reason,
$\tilde{\gamma}$ also has the properties $1)-3)$ like $\gamma_i$.

Choose $x=\tilde{\gamma}(a)\in \tilde{\gamma}\backslash\{p_0\}$ and
take $\gamma(s)=\tilde{\gamma}(a+s), s\geq 0$. We obtain a ray
$\gamma$ in $E$ emanating from  $x\in E$ as claimed.

Note for any $z\in \gamma$,  $d(z,D)\geq a>0$. So we may choose
small $R_0$   ($R_0<a$) such that $B_z(R_0)\subset E$, $z\in
\gamma$. Take $R_0$ satisfying $R_0< \min\{\frac{\pi}{2\sigma}, i_0,
a\}$. By Lemma \ref{volesti}, for any $z\in \gamma \subset M$,
\begin{equation}\label{volinequality}
\text{Vol}(B_z(R_0)\geq \sigma^{-n}\omega_n(\sin
R_0\sigma)^ne^{-H_0R_0}=\beta>0.
\end{equation}

Consider a sequence of points $z_j=\gamma (2jR_0), j=0,...,k-1,
\text{where}\quad k=\left[\frac{R}{2R_0}\right], R\geq 2R_0.$
Observe that any two balls $B_{z_j}(R_0)$ are disjoint and
$B_x(R)\supset \bigcup^{k-1}_{j=0}B_{z_j}(R_0)$. Then, $B_x(R)\cap
E\supset \bigcup^{k-1}_{j=0}B_{z_j}(R_0)$, and by
(\ref{volinequality}),
$$\text{Vol}(B_x(R)\cap E)\geq \text{Vol}(\bigcup^{k-1}_0B_{z_j}(R_0))\geq k\beta\geq(\frac{R}{2 R_0}-1)\beta, \quad R\geq 2R_0.$$
Hence $$\liminf_{R\rightarrow\infty}\frac{\text{Vol}(B_x(R)\cap
E)}{R}>0.$$

It is direct, from the definition of $\liminf$, that limit is
independent of the choice of $p$ and hence constant for any point of
$E$.

\end{proof}

\begin{cor}\label{volnonnegative} Let $N$ be a complete simply connected  manifold of nonpositive sectional curvature and $M$ be a complete noncompact  immersed submanifold in $N^m$
with norm-bounded mean curvature vector field $H$. Then each end of
$M$ has infinite volume.
\end{cor}

Li and Wang (\cite{lw}, Corollary 4) showed that if an end of a
manifold is of infinite volume and satisfies a Sobolev type
inequality, then this end must be nonparabolic. With this property,
we obtain Proposition \ref{nonparaend} and Proposition \ref{simnon}
as follows:

\begin{prop}\label{nonparaend}Let $N^{n+1}$ be a complete Riemannian manifold of bounded geometry and $M^n$ be a complete noncompact immersed CMC hypersurface
in $N$ with finite Morse index. If $\inf\tilde{\text{Ric}}>-nH^2$,
then each end of $M$ must be nonparabolic.
\end{prop}

\begin{proof} It is well known that a CMC hypersurface with finite Morse index is
strongly stable outside a compact domain (by the same argument in
\cite{fc}). Hence we assume that $M$ is stable outside a compact
domain $\Omega\subset M$. Clearly each end of $M$ is also stable
outside  $\Omega$.

 Since nonparabolicity of an end depends only on its infinity behavior, it is sufficient to show that each end $E$ of $M$ with respect to  any compact
 set $D$ $(\Omega\supset D)$ is nonparabolic.

  By stability,
 for any compactly supported function $f\in H_{1,2}(E)$, we have

 $$\int_E|\nabla f|^2\geq \int_E(\tilde{\text{Ric}}(\nu,\nu)+|\Phi|^2+nH^2)f^2\geq (\inf\tilde{\text{Ric}}+nH^2)\int_E f^2,$$

 that is, the end $E$ satisfies an Sobolev type inequality:
 $$\int_E f^2\leq C\int_E|\nabla f|^2.$$

By  Corollary 4 in \cite{lw} and Prop.\ref{infvol}, $E$ must be
nonparabolic.
\end{proof}

\begin{prop}\label{simnon}Let $N^m$ be a complete
 simply connected  manifold with nonpositive sectional curvature and let $M^n$ be a complete immersed minimal submanifold in $N^m$. If $n\geq 3$, then
each end of $M$ must be nonparabolic.
\end{prop}
{\bf Proof.} From  the theorem of Cartan-Hadamard, the exponential
map at any point of $N$ must be diffeomorphic $R^m$  and hence $N$
has bounded geometry. Assume $E$ is an end of $M$. Since under the
hypotheses of proposition, we have the following Sobolev inequality
(\cite{hs}, Theorem2.1):
\begin{equation}\label{sobo}(\int_E|f|^{\frac{2n}{n-2}})^{\frac{n-2}{n}}\leq C\int_E|\nabla
f|^2, f\in H_{1,2}(E).
\end{equation}

By Corollary 4 in \cite{lw} and Proposition \ref{infvol}, $E$ must
be nonparabolic. \qed
\begin{rema} The special case of Corollary \ref{volnonnegative} and Proposition \ref{simnon} that $E$
is an end of a minimal submanifold in $\mathbb{R}^m$ was proved in \cite{csz}.
\end{rema}

\section{Uniqueness of ends}

In this section we  discuss the uniqueness of ends of weakly stable
CMC hypersurfaces. We initially  prove an algebra inequality.

\begin{lem}
\label{lemalgebra}\label{lemma1}Let $A=(a_{ij})$ be an $n\times n$ real symmetric matrix with trace $%
\text{tr}(A)=nH$. Then
\begin{equation}
nHa_{11}-\sum_{i=1}^{n}a_{1i}^{2}\geq(n-1)H^{2}-(n-2)|H||B|\sqrt{\frac{n-1}{n}}-%
\frac{n-1}{n}|B|^{2},
\end{equation}
where $ B=(b_{ij})=A-HI, |B|^2=\sum_{i,j=1}^{n}b_{ij}{}^{2}$, $I$ is
the identity matrix.
\end{lem}
\begin{proof} Note $\sum_{i=1}^{n}b_{ii}=0$. We have
\begin{align*}
b_{11}^{2}  =\left( \sum_{i=2}^{n}b_{ii}\right) ^{2}
 \leq\left( n-1\right) \sum_{i=2}^{n}b_{ii}^{2}.
\end{align*}
Then
\begin{align}
|B|^{2} & =\sum_{i,j=1}^{n}b_{ij}{}^{2}  \geq
b_{11}^{2}+\sum_{i=2}^{n}b_{ii}^{2}+2\sum_{i=2}^{n}b_{1i}^{2} \notag
\\
& \geq b_{11}^{2}+\frac{1}{n-1}\left( \sum_{i=2}^{n}b_{ii}\right)
^{2}+2\sum_{i=2}^{n}b_{1i}^{2}  \notag \\
& \geq\frac{n}{n-1}\left(
b_{11}^{2}+\sum_{i=2}^{n}b_{1i}^{2}\right)\notag.
\end{align}

By $b_{ii}=a_{ii}-H, i=1,...,n$; $b_{ij}=a_{ij}$,  $i\neq j,
i,j=1,...,n,$  we have
\begin{align}\label{aeq}
nHa_{11}-\sum_{i=1}^{n}a_{1i}^{2}
& =(n-1)H^{2}+(n-2)Hb_{11}-(b_{11}^{2}+\sum_{i=2}^{n}b_{1i}^{2})\notag\\
&
\geq(n-1)H^{2}-(n-2)|H||B|\sqrt{\frac{n-1}{n}}-\frac{n-1}{n}|B|^{2}.
\end{align}

\end{proof}
 As a consequence, we obtain the following inequality,  which
was proved in \cite{chen} (Lemma 2.1 in \cite{chen}) by a different
proof.

\begin{prop}
\label{propalgebra}Let $A=(a_{ij})$ be an $n\times n$ real symmetric matrix with trace $%
\text{tr}(A)=nH$. Then
\begin{equation}\label{ab}
|A|^2+nHa_{11}-\sum_{i=1}^{n}a_{1i}^{2}\geq\frac{n^2(5-n)}{4}H^2,
\end{equation}
where $|A|^2=\sum_{i,j=1}^{n}a_{ij}^2$. Moreover  equality holds if
and only if one of the following cases occurs:

1) $n=2$, $A=HI$, where $I$ is the identity matrix;

2) $n\geq 3$, $A$ is a diagonal matrix with
$a_{11}=-\frac{n(n-1)}{2}H$, $a_{ii}=\frac{n}{2}H, i=2,...,n$, and
$a_{ij}=0, i\neq j, i,j=1,...,n$.

\end{prop}

\begin{proof} We use the same notations in Lemma \ref{lemma1}.
By $|B|^2=|A|^2-nH^2$,
\begin{align*}
&\quad |A|^2+nHa_{11}-\sum_{i=1}^{n}a_{1i}^{2} \\
& \geq |B|^2+(2n-1)H^{2}-(n-2)|H||B|\sqrt{\frac{n-1}{n}}-\frac{n-1}{n}|B|^{2}\\
&=(\frac{|B|}{\sqrt{n}}-\frac{(n-2)\sqrt{n-1}}{2}|H|)^2+\frac{n^2(5-n)H^2}{4}\\
&\geq \frac{n^2(5-n)H^2}{4}.
\end{align*}

Thus, we obtain (\ref{ab}). If the equality in (\ref{ab}) holds,
then,

1) in the case $n=2$,
$\frac{|B|}{\sqrt{n}}-\frac{(n-2)\sqrt{n-1}}{2}|H|=0$, so $B=0$,
that is,  $A=HI$.

2) in the case $n\geq 3$, by the proof of  (\ref{ab}), we have
$\sum_{i=1}^nb_{ii}=0$; $b_{ii}=b_{jj}, i,j=2,...,n$; $b_{ij}=0$,
$i\neq j,  i,j=1,...,n$;
$\frac{|B|}{\sqrt{n}}-\frac{(n-2)\sqrt{n-1}}{2}|H|=0$. Moreover,
$b_{11}$ and $H$  have different signs.

Thus,  $b_{11}=-(n-1)b_{22},
|b_{11}|=\frac{\sqrt{n-1}}{\sqrt{n}}|B|=\frac{(n-1)(n-2)}{2}|H|.$
Since $b_{11}$ and $H$ have  different signs,
$b_{11}=-\frac{(n-1)(n-2)}{2}H$.

Hence  $a_{11}=-\frac{n(n-1)}{2}H$, $a_{ii}=\frac{n}{2}H,
i=2,...,n$, and $a_{ij}=0,  i\neq j, i,j=1,...,n$, that is, $A$ is a
diagonal matrix with $a_{ii}$ given above.

Conversely, $A$ in (1) and (2) satisfy the equality in (\ref{ab}).
The proof is complete.

\end{proof}

Applying Proposition \ref{propalgebra} to  hypersurfaces, we obtain
the following Proposition \ref{secform}, which can be used to prove
Theorem \ref{thm3} and may have its independent interest.

\begin{prop}\label{secform} Let $N$ be an $(n+1)$-dimensional manifold
and $M$ be a hypersurface in $N$ with  mean curvature $H$ (not
necessarily constant). Then, for any local orthonormal frame
$\{e_i\}, i=1,...,n,$ of  $M$,
\begin{equation}\label{sf}
|A|^2+nHh_{11}-\sum^{n}_{i=1}h^2_{1i}\geq \frac{n^2(5-n)H^2}{4},
\end{equation}
where the second fundamental form $A=(h_{ij}),
h_{ij}=\left<Ae_i,e_j\right>, i,j=1,...,n$. Furthermore, the
equality in (\ref{sf}) holds for some $\{e_i\}$ at a point $p\in M$,
if and only if one of the following cases occurs at $p$:

(1) $n=2$, $A=HI$, where $I$ is the identity map, that is, $M$ is
umbilic at $p$;

(2) $n\geq 3$, $A$ is a diagonal matrix with
$a_{11}=-\frac{n(n-1)}{2}H$, $a_{ii}=\frac{n}{2}H, i=2,...,n$, and
$a_{ij}=0, i\neq j, i,j=1,...,n$, that is, $M$ has $n-1$ equal
principle curvature and only one is different when $H\neq 0$ at $p$,
or $M$ is totally geodesic when $H=0$ at $p$.
\end{prop}

Schoen and Yau (\cite{sy}, cf. \cite{lw})  proved an inequality on harmonic functions
on stable minimal hypersurfaces, we generalize their inequality
to the CMC $H$-hypersurfaces:

\begin{lem}\label{prop3}
Let $M$ be a complete  hypersurface with constant mean curvature $H$
in $N^{n+1}$. Suppose that $u$ is a harmonic function defined on
$M$. If $\varphi$ is a compactly supported  function $ \varphi\in
H_{1,2}(M)$ such that $\varphi |\nabla u|$ satisfies the stability
inequality (\ref{eqnstable}), then
\begin{align}\label{syineq}
\int_{M}\varphi^{2}|\nabla u|^{2}\{\dfrac{1}{n}&|\Phi|^{2}-\sqrt{\dfrac{n-1}{n}}(n-2)H|\Phi|+\left(2n-1\right)H^2\notag\\
&+\tilde{\mathrm{Ric}}\left(\frac{\nabla u}{\left|\nabla
u\right|},\frac{\nabla u}{\left|\nabla u\right|}\right)+
\tilde{\mathrm{Ric}}\left(\nu,\nu\right)-\tilde{\mathrm{K}}\left(\frac{\nabla
u}{\left|\nabla u\right|},\nu\right)\}
\\
&\quad\quad +\int_{M}\frac{1}{n-1}\varphi^{2}|\nabla|\nabla
u||^{2}\leq\int _{M}|\nabla\varphi|^{2}|\nabla u|^{2}.\notag
\end{align}

  \end{lem}
\begin{proof}Recall the Bochner formula
\begin{equation}\label{boch}\frac{1}{2}\Delta|\nabla
u|^{2}=\text{Ric}(\nabla u,\nabla u)+|\nabla
^{2}u|^{2},
\end{equation}
\begin{equation}\label{equaharm}\textrm{the equality}\quad\qquad\quad \frac{1}{2}\Delta|\nabla u|^{2}=|\nabla u|\Delta|\nabla u|+|\nabla|\nabla u||^{2},
\end{equation}
and the inequality (see \cite{lw}): when $u$ is harmonic function,
\begin{equation}\label{ineqharm}
|\nabla^{2}u|^{2}\geq\frac{n}{n-1}|\nabla|\nabla u||^{2}.
\end{equation}
By (\ref{boch}), (\ref{equaharm}) and (\ref{ineqharm}), we have
$$|\nabla u|\Delta|\nabla u|\geq \mathrm{Ric}(\nabla u,\nabla u)+\frac{1}{n-1}|\nabla|\nabla u||^{2}.$$

Let $\varphi$ be a locally Lipschitz function with compact support
on $M$. Choose $f=\varphi |\nabla u|$ in the stability inequality
(\ref{eqnstable}). Then
\begin{align}\label{ieq}
&\int_{M}( |\Phi|^{2}+\tilde{\mathrm{Ric}}\left( \nu,\nu\right)
+nH^{2}) \varphi^{2}|\nabla u|^{2}\notag \\
 \leq&\int_{M}|\nabla\left(
\varphi\left| \nabla u\right| \right) |^{2} \notag\\
=&\int_{M}|\nabla\varphi|^{2}|\nabla u|^{2} -2\langle \varphi\left(
\nabla\left| \nabla u\right| \right) ,\left| \nabla u\right|
\nabla\varphi\rangle +\int\varphi^{2}\left|
\nabla|\nabla u|\right| ^{2} \notag\\
= & \int_{M}|\nabla\varphi|^{2}|\nabla u|^{2}-\int\varphi^{2}|\nabla
u|\Delta|\nabla u| \notag\\
 \leq&\int_{M}|\nabla\varphi|^{2}|\nabla u|^{2}-\frac{1}{n-1}\int_M\varphi^2|\nabla|\nabla u||^{2}-\int_M\varphi^2\mathrm{Ric}(\nabla u,\nabla u).
\end{align}

For any point $p\in M$ and any unit vector $\eta\in T_{p}M$, we
choose a local orthonormal frame field $\{e_{1},e_{2},\cdots,e_{n}\}$ \ at $%
p $ such that $e_{1}=\eta$, we have, from Gauss equation:
\begin{equation}
K(e_{i},e_{j})-\tilde{K}(e_{i},e_{j})=h_{ii}h_{jj}-h_{ij}^{2},\text{ for }%
i,j=1,2,\cdots ,n,  \label{eqngauss}
\end{equation}
\begin{align}\label{eqric}
{\mathrm{Ric}}(\eta,\eta) & =\sum_{i=2}^{n}\tilde{K}(\eta,e_{i})+h_{11}%
\sum_{i=2}^{n}h_{ii}-\sum_{i=2}^{n}h_{1j}^{2} \notag  \\
& =\tilde{ \mathrm{Ric}}\left(\eta,\eta\right) -\tilde{K}\left(
\nu,\eta\right)+nHh_{11}-\sum_{i=1}^{n}h_{1j}^{2}.
\end{align}

Substituting  $\eta=\frac{\nabla u}{|\nabla u|}$ into (\ref{eqric})
and then substituting (\ref{eqric}) into  (\ref{ieq}), we obtain
\begin{align}\label{beq}
&\int_{M}( |\Phi|^{2}+\tilde{\mathrm{Ric}}\left( \nu,\nu\right)
+nH^{2}) \varphi^{2}|\nabla u|^{2}\notag \\
 \leq&\int_{M}|\nabla\varphi|^{2}|\nabla u|^{2}-\int_{M}\frac{1}{n-1}\varphi^{2}|\nabla|\nabla
 u||^{2}\notag
\\
 \quad &-\int_M\varphi^2\{\tilde{\mathrm{Ric}}\left( \frac{\nabla u}{\left| \nabla u\right|},\frac{\nabla u}{\left| \nabla u\right|}\right) -
\tilde{K}\left( \nu,\frac{\nabla u}{\left| \nabla
u\right|}\right)+nHh_{11}-\sum_{i=1}^{n}h_{1i}^{2} \}.
\end{align}

By Lemma \ref{lemma1},
\begin{align*}
&\int_{M}( |\Phi|^{2}+\tilde{\mathrm{Ric}}\left( \nu,\nu\right)
+nH^{2}) \varphi^{2}|\nabla u|^{2}\notag \\
  \leq&\int_{M}|\nabla\varphi|^{2}|\nabla u|^{2}-\int_{M}\frac{1}{n-1}\varphi^{2}|\nabla|\nabla u||^{2}
\\
& +\int_{M}\varphi
^{2}|\nabla u|^{2}\{(n-2)H|\Phi|\sqrt{\frac{n-1}{n}}+\frac{n-1}{n}|\Phi|^{2}-(n-1)H^2\\
&-\tilde{\mathrm{Ric}}\left( \frac{\nabla u}{\left| \nabla
u\right|},\frac{\nabla u}{\left| \nabla u\right|}\right) +
\tilde{K}\left( \nu,\frac{\nabla u}{\left| \nabla
u\right|}\right)\}.
\end{align*}
Thus (\ref{syineq}) holds.
\end{proof}

As  mentioned in Section \ref{pre}, in order to estimate the number
of the ends of  hypersurface $M$, we need to discuss the
nonexistence of nonconstant bounded harmonic functions on $M$ with
finite Dirichlet integral. We obtain that
\begin{thm}\label{thm3}
Let $M$ be a complete noncompact weakly stable CMC $H$-hypersurface in $N^{n+1}$ in a manifold $N$. If, for any $p\in M$,
$$\tilde{\mathrm{Ric}}\left(X,X\right)+\tilde{\mathrm{Ric}}\left(\nu,\nu\right)-\tilde{K}\left(X,\nu\right)\geq \frac{n^2(n-5)}4H^2, X\in T_pM, |X|=1,$$
then $M$ does not admit nonconstant bounded harmonic functions with finite Dirichlet integral.

\end{thm}

\begin{proof} We prove the conclusion by contradiction.  Suppose there exists a nonconstant bounded harmonic function $u$ with finite Dirichlet
integral on $M$. Then there exists some  point $p\in M$ such that
$|\nabla u|(p)\neq 0$. Hence, $\int_{B_p(a)}|\nabla u|>0$, for all
$a>0$.

{\it We claim that $u$ must satisfy $\int_M|\nabla u|=\infty$.}

 By the boundness of $u$,
$\int_{B_p(R)}|\nabla u|^2=\int_{\partial B_p(R)}u\frac{\partial
u}{\partial r} \leq C\int_{\partial B_p(R)}|\nabla u|,$ where $C$ is
a constant. Hence when $R>1$, $$0<C_0=\int_{B_p(1)}|\nabla
u|^2\leq\int_{B_p(R)}|\nabla u|^2\leq C\int_{\partial B_p(R)}|\nabla
u|,$$

that is $\int_{\partial B_p(R)}|\nabla u|\geq C_1>0.$

By co-area formula,
\begin{equation}
\int_{B_p(R)}|\nabla u|=\int_1^Rdr\int_{\partial B_p(r)}|\nabla
u|\geq C_1(R-1).
\end{equation}
Letting $R\rightarrow\infty$, we have $\int_M|\nabla u|=\infty$ as
claimed.

Take, for $R>a$,
\begin{equation}
\varphi_1(a,R)=
\begin{cases}
1,&\text{ on }\bar{B}_{p}(a), \\
\frac{a+R-x}{R},&\text{ on }B_{p}(a+R)\backslash B_{p}(a), \\
0,&\text{ on }M\backslash B_{p}(a+R).%
\end{cases}
\end{equation}
and
\begin{equation}
\varphi_2(a,R)=
\begin{cases}
0,&\text{ on }B_{p}(a+R), \\
\frac{a+R-x}{R},&\text{ on }B_{p}(a+2R)\backslash B_{p}(a+R), \\
-1,&\text{ on }B_p(a+2R+b)\backslash B_{p}(a+2R),\\
\frac{x-(a+3R+b)}{R},&\text{ on }B_p(a+3R+b)\backslash B_{p}(a+2R+b),\\
0,&\text{ on }M\backslash B_{p}(a+3R+b),
\end{cases}
\end{equation}
where constant $b>0$ will be determined later.

For any $\epsilon>0$ given, we may choose large $R$ such that $\frac{1}{R^2}\int_M|\nabla u|^2<\epsilon$.

Define $\psi(t,a,R)=\varphi_1(a,R)+t\varphi_2(a,R)$, $t\in [0,1]$.
We have
$$\int_M\psi(0,a,R)|\nabla u|\geq\int_{B_p(a)}|\nabla u|>0,$$ and
\begin{align}
\int_M\psi(1,a,R)|\nabla u|&=\int_M(\varphi_1(a,R)+\varphi_2(a,R))|\nabla u|\notag\\
&\leq\int_{B_p(a+R)}|\nabla u|-\int_{B_p(a+2R+b)\backslash B_p(a+2R)}|\nabla u|.
\end{align}

By claim, for $a$ and $R$ fixed, we may find $b$ sufficiently large, depending on $a$ and $R$ such that
$$\int_M\psi(1,a,R)|\nabla u|<0.$$
By the continuity of $\psi(t,a,R)$ on $t$, there exists some $t_0\in
(0,1)$ depending on $a$ and $R$ such that
$\int_M\psi(t_0,a,R)|\nabla u|=0.$

Denote $\psi(t_0,a,R)$ by $f$. Since $M$ is weakly stable,
$f=\psi(t_0,a,R)|\nabla u|$ satisfies the stability  inequality
(\ref{eqnstabletr}) and   hence also satisfies  Lemma \ref{prop3}.

 Note the
curvature condition  implies that
\begin{align*}
\dfrac{1}{n}|\Phi|^{2}&-\sqrt{\dfrac{n-1}{n}}(n-2)H|\Phi|+\left(2n-1\right) H^{2}\\
&+\tilde{\mathrm{Ric}}\left(\frac{\nabla u}{\left|\nabla u\right|},\frac{\nabla u}{\left|\nabla u\right|}\right)+
\tilde{\mathrm{Ric}}\left(\nu,\nu\right)-\tilde{\mathrm{K}}\left(\frac{\nabla u}{\left|\nabla u\right|},\nu\right)\\
&\quad\geq\dfrac{1}{n}|\Phi|^{2}-\sqrt{\dfrac{n-1}{n}}(n-2)H|\Phi|+\frac{(n-1)(n-2)^2}{4}H^2\\
&\quad\geq[\frac{|\Phi|}{\sqrt{n}}-\frac{\sqrt{n-1}(n-2)H}{2}]^2\geq 0.
\end{align*}

Then by (\ref{syineq}), we have
\begin{equation}
\int_{M}\frac{1}{n-1} f^2|\nabla|\nabla u||^{2}\leq\int _{M}|\nabla
f|^{2}|\nabla u|^{2}.
\end{equation}

 Then
\begin{align*}
&\frac{1}{n-1}\int_{B_p(a)}|\nabla|\nabla
u||^2\notag\\&\leq\int_{B_{p}(a+2R)\backslash
B_{p}(a)}|\nabla\varphi_1|^2|\nabla u|^2+t_0^2
\int_{B_p(a+3R+b)\backslash B_{p}(a+2R+b)}|\nabla\varphi_2|^2|\nabla u|^2\\
&\leq\frac{1}{R^2}\int_{B_{p}(a+2R)\backslash B_{p}(a)}|\nabla u|^2+\frac{1}{R^2}\int_{B_p(a+3R+b)\backslash B_{p}(a+2R+b)}|\nabla u|^2\\
&\leq\frac{1}{R^2}\int_M|\nabla u|^2<\epsilon.
\end{align*}
In the above first inequality, we used
$\left\langle\nabla\varphi_{1},\nabla\varphi_{2}\right\rangle=0$.

By the arbitrariness of $\epsilon$ and $a$, $\nabla|\nabla u|\equiv
0$. So $|\nabla u|\equiv constant$.

If $|\nabla u|\equiv const.\neq 0$, then $u$ is a nonconstant
bounded harmonic function. This says $M$ must be nonparabolic. Thus
$\mathrm{vol}(M)=\infty$. Hence $\int_M|\nabla u|^2=\infty$, which
is impossible. Therefore $|\nabla u|\equiv  0$, $ u\equiv constant$.
Contradiction. The proof is complete.
\end{proof}
\begin{rema} If $N$ is $3$-dimensional, the curvature
$\tilde{\mathrm{Ric}}(X)+\tilde{\mathrm{Ric}}(Y)-\tilde{\mathrm{K}}(X,Y),
 X, Y\in T_pN,  X\perp Y, |X|=|Y|=1, p\in N,$ is equal to  the scalar
curvature $\tilde{\mathrm{S}}$. From the definition we know that the
nonnegativity of the sectional curvature  implies the nonnegativity
of the above curvature. However, there are some examples showing
that the converse may not be true (see \cite{shy}). In this paper,
we adopt this curvature because it appears  naturally in this
context and provides more examples.
\end{rema}

Now  we are ready to obtain the uniqueness of the ends of weakly
stable CMC hypersurfaces. First, we consider the weakly stable
minimal hypersurfaces in $\mathbb{R}^{n+1}$ and obtain

\begin{thm} \label{oneRRn}(Th.\ref{oneRn1}) If $M$ is a complete oriented weakly stable minimal
immersed hypersurface in $\mathbb{R}^{n+1}$, $n\geq 3$,
 then $M$ must have only one end.
\end{thm}
\begin{proof} First a complete minimal hypersurface in $\mathbb{R}^{n+1}$ must be compact.
By Theorem \ref{thm3}, the dimension of the space
$\mathcal{H}^0_D(M)$ is $1$. By Proposition \ref{simnon}, each end
of $M$ must be nonparabolic. Hence by Theorem \ref{endform}, $M$
must have only one end.
\end{proof}

Recall that Anderson (\cite{an}, Theorem 5.2) proved that a complete
minimal hypersurface in $\mathbb{R}^{n+1}  (n\geq 3)$ with finite
total curvature and one end must be an affine-plane. Hence by the
result of Anderson and Theorem \ref{oneRRn}, we have
\begin{cor} \label{plane}(Cor.\ref{plane1}) A complete weakly stable immersed minimal hypersurface $M$ in $\mathbb{R}^{n+1}, n\geq 3,$ with
finite total  curvature (i.e., $\int_M|A|^n<\infty$) is a
hyperplane.
\end{cor}
\begin{rema} Y.B. Shen and X. Zhu \cite{shz} showed that a complete  stable immersed minimal hypersurface in $\mathbb{R}^{n+1}$ with finite total
 curvature  is a hyperplane. So Corollary \ref{plane} generalizes
their result.
\end{rema}
\begin{cor} A complete weakly stable CMC hypersurface in
$\mathbb{R}^{n+1}, n\geq 3,$ with finite total curvature (i.e.,
$\int_M|\Phi|^n<\infty$) is either a hyperplane or a geodesic
sphere.
\end{cor}
\begin{proof} do Carmo, Cheung
and Santos \cite{dccs} proved that a complete stable CMC
$H$-hypersurface, $H\neq 0$, in $\mathbb{R}^{n+1}, n\geq 3$ with
finite total curvature  must be compact. By their result, Theorem2.1
in \cite{bc} and Corollary \ref{plane}, we obtain that a complete
weakly stable CMC hypersurfaces in $\mathbb{R}^{n+1}$ with finite
total curvature must be a hyperplane or a geodesic sphere.
\end{proof}

When the Ricci curvature of the ambient manifold has a strict low
bound $-nH^2$, we obtain, for CMC $H$-hypersurfaces,

\begin{thm}\label{one} Let $N^{n+1}$ be a complete Riemannian manifold of bounded geometry and $M$ be a complete noncompact weakly stable  immersed
CMC $H$-hypersurface in $N$. If $\inf\tilde{\text{Ric}}>-nH^2$ and
for any $p\in M$, $X\in T_pM, |X|=1$,
\begin{equation}\label{bric}
\tilde{\mathrm{Ric}}\left(X,X\right)+\tilde{\mathrm{Ric}}\left(\nu,\nu\right)-\tilde{K}\left(X,\nu\right)\geq
\frac{n^2(n-5)}4H^2,
\end{equation}
then $M$ has only one end.
\end{thm}
\begin{proof} By Proposition \ref{nonparaend}, each end of $M$ is nonparabolic. Hence, by Theorem \ref{thm3} and  Theorem \ref{endform}, we get the conclusion.
\end{proof}
\begin{rema}The curvature conditions in Theorem \ref{one} demands,
 some positivity of curvature of $N$ or the restriction on the
 dimension of $N$. See the evidence in its several consequences. However,
 it is worth to note that  the same result holds
 without any curvature condition for minimal
 hypersurfaces in $\mathbb{R}^{n+1}, n\geq 3$, as Theorem \ref{oneRRn} says. The reason
 is that we have a global Sobolev inequality (\ref{sobo}) in this
 case.

\end{rema}

 Now we derive some consequences of Theorem \ref{one}.
\begin{cor}\label{oneHn}(Th.\ref{oneHn1}) Let $M$ be a complete noncompact weakly stable CMC $H$-hypersurface in the hyperbolic space
$\mathbb{H}^{n+1}, n=3,4$. If $H^2\geq\frac{10}{9}$, when $n=3$;
$H^2\geq\frac{7}{4}$, when $ n=4$, respectively, then $M$ has only
one end.
\end{cor}

\begin{proof} Observe that the hypotheses of  Theorem \ref{one} are satisfied.
\end{proof}

 In \cite{ch1} and
 \cite{ch2}, Cheng
 proved   that let $N^{n+1}$ be an $(n+1)$-dimensional manifold  and $M$ be a complete CMC $H$-hypersurface
immersed in $N$ with finite index. Then $M$ must be compact in the
following two cases: 1) $n=3,4,5$, $H=0$,
$\inf\{\tilde{\mathrm{Ric}}(w)+\tilde{\mathrm{Ric}}(\nu)-\tilde{\mathrm{K}}(w,\nu)|
\,\,w\in T_p^1M,\,  p\in M\}> 0$; 2) $n=3,4$, $H\neq 0$,
$\inf\{\tilde{\mathrm{Ric}}(w)+\tilde{\mathrm{Ric}}(\nu)-\tilde{\mathrm{K}}(w,\nu)|
\,\,w\in T_p^1M,\,  p\in M\}> \frac{n^2(n-5)}{4}H^2$. Combining this
result, Theorem \ref{oneRRn} and Theorem \ref{one}, we obtain that

\begin{cor}\label{oneRn} If $M$ is a complete noncompact weakly stable
immersed CMC $H$-hypersurface in $\mathbb{R}^6$,
 then $M$ has only one end.
\end{cor}
\begin{proof} If $H\neq 0$,  by Theorem \ref{one}, we know that when $n\leq 5$, $M$ in $\mathbb{R}^{n+1}$ has only one end.  But it is
known (\cite{ch1}) that there is no complete noncompact weakly
stable CMC $H$-hypersurfaces in $\mathbb{R}^4, \mathbb{R}^5$ ($H\neq
0$). Hence only the case $\mathbb{R}^6$ may occur. If $H=0$,
 Theorem \ref{oneRRn}  says that a complete noncompact weakly stable minimal hypersurfaces in $\mathbb{R}^{n+1}, n\geq 3,$ must have only one end.
 Combining two cases, we obtain  the conclusion.
\end{proof}

The above result in \cite{ch1} and \cite{ch2}  implies that any
complete weakly stable  $H$-hypersurface in a complete
   manifold $N^{n+1}$, $n=3,4$, or $n=5$ and $H=0$,  of nonnegative sectional curvature must be
   compact. Hence by Theorem \ref{one}, we obtain the following

\begin{thm}\label{sectcurvature} (Th.\ref{sectcurvature1}) Let $N^{n+1}, n\geq 5,$ be a complete Riemannian manifold and $M$ be a complete
noncompact weakly stable immersed  CMC $H$-hypersurface in $N$. If
one of the following cases  occurs,

(1) when $n=5$, the sectional curvature of $N$ is nonnegative and $H\neq 0$;

(2) when $n\geq 6$, the sectional curvature $\tilde{K}$ of $N$
satisfies $\tilde{K}\geq \tau>0$, and $H^2\leq
\frac{4(2n-1)}{n^2(n-5)}\tau$ for some  number $\tau>0$;

 (3) when $n\geq 6$, the sectional curvature and the Ricci curvature of
$N$ satisfy $\tilde{K}\geq 0$, $\tilde{\mathrm{Ric}}\geq \tau
>0$,  and $H^2\leq \frac{4\tau}{n^2(n-5)}$, for some number $\tau
>0$,

 then  $M$ has only one end.

 In particular, any complete noncompact
 stable minimal hypersurface in a manifold $N^{n+1}$, $n\geq 6,$ of nonnegative sectional
 curvature and Ricci curvature bounded from below by a positive
 number has only one end.

\end{thm}

As some special cases, Theorem \ref{sectcurvature} implies that

\begin{cor}\label{oneSn} A complete noncompact weakly stable CMC
$H$-hypersurface has only one end, if it is in either

1) the standard sphere $\mathbb{S}^6$ with $H\neq 0$; or

2)  the standard sphere $\mathbb{S}^{n+1}, n\geq 6,$ with $H^2\leq
\frac{4(2n-1)}{n^2(n-5)}$; or

3) $\mathbb{S}^k\times \mathbb{S}^l,$  $k\geq 2, l\geq2,$  $k+l\geq
7,$ with the product metric and $H=0$.

In particular, a complete noncompact stable minimal hypersurface in
$\mathbb{S}^{n+1},$  $n\geq 6,$ and $\mathbb{S}^k\times
\mathbb{S}^l,$ $k\geq 2, l\geq2, k+l\geq 7,$ has only one end.
\end{cor}

\section{Property of Schr\"odinger operator on parabolic
manifolds}\label{sectschrodinger}

In this section, we   prove a property of Schr\"odinger operator for
parabolic manifolds (not necessary a submanifold). It will be
applied to weakly stable CMC hypersurfaces and also may have its
independent interest.
\begin{prop}\label{sch}
Let $M$ be a complete  parabolic   manifold with infinity volume.
Consider the operator $L=\Delta+q(x)$ on $M$ (here $q:M\rightarrow
\mathbb{R}$ is a continuous function on $M$).  If $q(x)\ge 0 $ and
$q(x)$ is not identically zero, then there exists a compactly
supported piecewise smooth function $\psi$ such that
$\int_M\psi(x)=0$ and $-\int_M\psi L\psi<0$.
\end{prop}
\begin{proof} By  hypothesis, we may choose a point $p\in M$ such that $q(p)>0$. Denote $C:=\int_{B_p(1)}q(x) dv>0$.
 Choose a monotonically increasing sequence $\{r_i\}$ with
$r_i\rightarrow \infty$ and consider  the harmonic functions $g_i$
defined by
$$
\begin{cases}
\Delta g_{i}  =0,&\text{ on }B_{p}(r_{i})\backslash B_{p}(1), \\
g_{i} =1,&\text{ on }\partial B_{p}(1), \\
g_{i} =0,&\text{ on }\partial B_{p}(r_{i}).
\end{cases}
$$
Since $M$ is parabolic, we have that $\lim_{r_i\to
+\infty}\int_{B_{p}(r_{i})\backslash B_{p}(1)}|\nabla g_i|^2=0.$

By this property, we can find some positive number $R_1>1$  and a
corresponding function $f_1$ satisfying
$$
\begin{cases}
\Delta f_{1}  =0,&\text{ on }B_{p}(R_{1})\backslash B_{p}(1), \\
f_{1} =1,&\text{ on }\partial B_{p}(1), \\
f_{1} =0,&\text{ on }\partial B_{p}(R_{1}),
\end{cases}
$$
and
 $\int_{B_{p}(R_{1})\backslash B_{p}(1)}|\nabla f_1|^2<\frac{C}6.$

  Let
\begin{equation}
\varphi_1=
\begin{cases}
1,&\text{ on }\bar{B}_{p}(1), \\
f_{1},&\text{ on }B_{p}(R_{1})\backslash B_{p}(1), \\
0,&\text{ on }M\backslash B_{p}(R_{1}).%
\end{cases}
\end{equation}

Similarly we can find a positive number $R_2>R_1$ and a  function
$f_2$ satisfying $\int_{B_{p}(R_2)\backslash B_{p}(R_1)}|\nabla
f_2|^2<\frac{C}6$ and
$$
\begin{cases}
\Delta f_{2}  =0,&\text{ on }B_{p}(R_{2})\backslash B_{p}(R_1), \\
f_{2}  =1,&\text{ on }\partial B_{p}(R_1), \\
f_{2} =0, &\text{ on }\partial B_{p}(R_{2}).
\end{cases}
$$
Let

\begin{equation}
\varphi_2=
\begin{cases}
0,&\text{ on }\bar{B}_{p}(R_1), \\
f_{2}-1,&\text{ on }B_{p}(R_2)\backslash B_{p}(R_1), \\
-1,&\text{ on }M\backslash B_{p}(R_2).%
\end{cases}
\end{equation}

Again, for any constant $b>0$,  there exists $R_3>R_2+b$ and a
function $f_3$ satisfying  $\int_{B_{p}(R_{3})\backslash
B_{p}(R_2+b)}|\nabla f_3|^2<\frac{C}6$ and

\begin{equation*}
\begin{cases}
\Delta f_{3}  =0,&\text{ on }B_{p}(R_{3})\backslash B_{p}(R_2+b), \\
f_{3}  =-1,&\text{ on }\partial B_{p}(R_2+b), \\
f_{3}  =0,&\text{ on }\partial B_{p}(R_{3}).
\end{cases}
\end{equation*}
Let
\begin{equation}
\varphi_3=
\begin{cases}
0,&\text{ on }\bar B_{p}(R_2+b), \\
f_{3}+1,&\text{ on }B_{p}(R_3)\backslash B_{p}(R_2+b), \\
1,&\text{ on }M\backslash B_{p}(R_3).%
\end{cases}
\end{equation}

Thus the sum of two functions $ \varphi_2+\varphi_3$ satisfies
\begin{equation}
\varphi_2+\varphi_3=
\begin{cases}
0,&\text{ on }B_{p}(R_1), \\
\varphi_{2},&\text{ on }B_{p}(R_2)\backslash B_{p}(R_1), \\
-1,&\text{ on }B_p(R_2+b)\backslash B_{p}(R_2),\\
f_3,&\text{ on }B_p(R_3)\backslash B_{p}(R_2+b),\\
0,&\text{ on }M\backslash B_{p}(R_3).
\end{cases}
\end{equation}
Let $\phi_t=\varphi_1+t( \varphi_2+\varphi_3)$. We see that $\phi_t
$ has compact support in $M$. Then we define
$\xi(t):=\int_M\varphi_1+t\int_M (\varphi_2+\varphi_3)$ we know that
$\xi(0)=\int_M\varphi_1>0 $ and since the volume of $M$ is infinite
we can choose $b$ large  such that
\begin{align*}
\xi(1)&=\int_M\varphi_1+\int_M( \varphi_2+\varphi_3)\\
&\leq \int_M\varphi_1-\int_{B_p(R_2+b)\backslash B_{p}(R_2)}1<0.
\end{align*}

So there exists a $t_0\in (0,1)$ such that $\int_M\phi_{t_0}=0$ and
\begin{align*}
-&\int_M\phi_{t_0} L\phi_{t_0}
=\int_M|\nabla \phi_{t_0}|^2-q(x)\phi_{t_0}^2\\
\le&\int_M|\nabla \varphi_{1}|^2+\int_M|\nabla \varphi_{2}|^2+\int_M|\nabla \varphi_{3}|^2-\int_{B_p(1)}q(x) \varphi_{1}^2\\
=&\int_{B_{p}(R_1)\backslash B_{p}(1)}|\nabla f_{1}|^2+\int_{B_{p}(R_2)\backslash B_{p}(R_1)}|\nabla f_{2}|^2\\
&\qquad\qquad\qquad+\int_{B_p(R_3)\backslash B_{p}(R_2+b)}|\nabla f_{3}|^2-\int_{B_p(1)}q(x)\\
<&\frac{C}6+\frac{C}6+\frac{C}6-C=-\frac{C}2<0.
\end{align*}
Choosing $\phi_{t_0}$ as $\psi$, we obtain the conclusion of the
proposition.
\end{proof}

\begin{rema} A special case of Proposition \ref{sch} that $M$ is a surface was proved by  da Silveira \cite{ds} by using  the
conformal structure of ends of two-dimensional parabolic manifolds.
This structure (obtained by using Huber's theorem)  does not exist
in higher dimensional cases.
\end{rema}

\section{Structure of weakly stable CMC hypersurfaces}\label{sectsch}

 In this section, we will study the structure of a weakly stable
CMC hypersurface according to its parabolicity or nonparabolicity.

{\bf (I) Parabolic case:}

Applying Proposition \ref{sch} to the case that $M$ is a weakly
stable CMC hypersurface, we obtain
\begin{prop}\label{prop2}
Let $M$ be a complete weakly stable CMC $H$-hypersurface in $N^{n+1}$. Suppose that the Ricci curvature of $N$ is bounded
from below by $-nH^{2}$. If $M$ is parabolic and has infinite volume, then $M$ must be totally umbilic in
$N.$ Moreover the Ricci curvature Ric$(\nu,\nu)$  in the normal direction is
identically equal to $-nH^{2}$ along $M$ and the scalar curvature $S_{M}$ is nonnegative.
\end{prop}

\begin{proof} From the assumption, $|\Phi|^2+\tilde{\mathrm{Ric}}(\nu ,\nu )+nH^{2}\ge 0$. Since $M$ is weakly stable, by Proposition \ref{sch}, it holds
that $|\Phi|^2+\tilde{\mathrm{Ric}}(\nu ,\nu )+nH^{2}\equiv 0$. Hence $\Phi\equiv 0$, that is, $M$ is umbilic, and $\tilde{\mathrm{Ric}}(\nu ,\nu )+
nH^{2}\equiv 0$.

At any point $p\in M$, choose a local orthonormal frame field $
e_{1},e_{2},\cdots ,e_{n},\nu $ \ at $p$ such that
$e_{1},e_{2},\cdots ,e_{n}$ are tangent fields.

Since $\Phi \equiv 0$, Gauss equations (\ref{eqngauss}) become:
\begin{equation}\label{gauss}
K(e_{i},e_{j})-\tilde{K}(e_{i},e_{j})=H^{2},\text{ when }i\neq j.
\end{equation}%
Then
\begin{equation}
\sum_{i,j=1}^{n}K(e_{i},e_{j})-\sum_{i,j=1}^{n}\tilde{K}%
(e_{i},e_{j})-n(n-1)H^{2}=0,
\end{equation}

\begin{align*}
S_{M} & =\sum_{i=1}^{n}[\tilde{\mathrm{Ric}}(e_{i},e_{i})-\tilde {K}%
(\nu,e_{i})]+n(n-1)H^{2} \\
& \geq-n^{2}H^{2}-\tilde{\mathrm{Ric}}(\nu,\nu)+n(n-1)H^{2} \\
& =0.
\end{align*}
\end{proof}

\begin{thm}\label{thm21}
Let $M$ be a complete weakly stable CMC hypersurface immersed in
$N^{n+1}$ with constant mean curvature $H$. Suppose that $N$ has
bounded geometry and the Ricci curvature of $N$ is bounded from
below by $-nH^{2}$. If $M$ is parabolic, then $M$ must be totally
umbilic in $N.$ Moreover the Ricci curvature Ric$(\nu,\nu)$ \ in the
normal direction is identically equal to $-nH^{2}$ along $M$ and the
scalar curvature $S_{M}$ is nonnegative.
\end{thm}

\begin{proof}
 From Lemma \ref{infvol} when $N$ has bounded geometry, then the volume of $M$ is infinite. Thus the conclusion follows directly from Proposition
  \ref{prop2}.
\end{proof}

\begin{thm}\label{thmpara}(Th.\ref{thmpara1}) Let $N$ be a complete  manifold of bounded geometry and $M$ be
a complete noncompact  weakly stable CMC $H$-hypersurface immersed
in $N$. If the sectional curvature of $N$ is bounded from below by
$-H^2$ and if $M$ is parabolic, then it is totally umbilic and has
nonnegative sectional curvature. Further,  either

(1) $M$ has only one end; or

(2) $M=\mathbb{R}\times P$ with the product metric, where $P$ is a
compact manifold of nonnegative sectional curvature.

\end{thm}

\begin{proof} We have shown $\Phi\equiv 0$ in Theorem \ref{thm21}. Since $\tilde{K}\geq -H^2$, $M$ has nonnegative (intrinsic) sectional curvature by the Gauss
 equation (\ref{gauss}). If $M$ has more than one end, by the splitting theorem of Cheeger and Gromoll \cite{cg} on manifolds of nonnegative curvature,
  we get the conclusion (2).
\end{proof}

{\bf (II) Nonparabolic case:}

In this situation, we  apply Theorem \ref{bric} and obtain the
following result:

\begin{thm}\label{thmnonpara}(Th.\ref{thmnonpara1}) Let $N$ be a complete Riemannian manifold and $M$ be
a complete noncompact  weakly stable $H$-hypersurface immersed in $N$. If $M$ is nonparabolic, and
$$\tilde{\mathrm{Ric}}(\nu)+\tilde{\mathrm{Ric}}(X)-\tilde{\mathrm{K}}(X,\nu)\geq \frac{n^2(n-5)}{4}H^2, \forall X\in T_pM, |X|=1, p\in M,$$
then it has only one nonparabolic end.
\end{thm}

\begin{proof} Since $M$ is nonparabolic, it has at least a nonparabolic end. If $M$ has  two or more nonparabolic ends,
then the dimension of $\mathcal{H}_D^0 (M)$ is not less than $2$,
which is a contradiction with Theorem \ref{thm3}.

\end{proof}

When $M$ is a weakly stable minimal hypersurface, combining Theorem
\ref{thmpara} in (I) and Theorem \ref{thmnonpara} in (II), we obtain
that,
\begin{thm}\label{structure} Let $N$ be a complete Riemannian manifold of bounded geometry and nonnegative sectional curvature and $M$ be
a complete noncompact oriented weakly stable minimal hypersurface
immersed in $N$. Then

(1) when $M$ is parabolic, then  either it has only one end and
nonnegative curvature; or it is isometric to $\mathbb{R}\times P$
with the product metric, where $P$ is a compact manifold of
nonnegative curvature. Moreover $M$ is totally geodesic;

(2) when $M$ is nonparabolic, then it has only one nonparabolic end.

\end{thm}

\bigskip\bigskip\noindent Xu Cheng\newline
Insitituto de Matematica\newline
Universidade Federal Fluminense-UFF\newline
Centro, Niter\'oi, RJ 24020-140\newline
Brazil\newline
email: xcheng@impa.br\\

\bigskip\bigskip\noindent Leung-fu Cheung\newline Department of
Mathematics\newline The Chinese University of Hong Kong\newline Shatin, Hong Kong\newline e-mail: lfcheung@math.cuhk.edu.hk\\

\noindent Detang Zhou\newline Insitituto de Matematica\newline
Universidade Federal Fluminense- UFF\newline Centro, Niter\'{o}i, RJ
24020-140, Brazil \newline email: zhou@impa.br
\smallskip
\newline{\it Current address:}\newline
Department of Mathematics\newline University of California,
Irvine\newline Irvine, CA 92697\newline USA

\enddocument

\end{document}